\documentclass[11pt]{amsart}

\usepackage{amscd}
\usepackage{amsmath}
\usepackage{graphicx}
\usepackage{amsfonts}
\usepackage{amssymb}
\begin{document}

\title{STRONG $3$-COMMUTATIVITY PRESERVING MAPS ON STANDARD OPERATOR ALGEBRAS}

\author{Meiyun Liu, Jinchuan Hou}
\address{Department of Mathematics, Taiyuan University of Technology ,
Taiyuan 030024, P. R. of China} \email{liumeiyunmath@163.com;
jinchuanhou@aliyun.com}


\thanks{This work is partially supported by  Natural
Science Foundation of China (11171249, 11271217).}

\begin{abstract}

Let $X$ be a Banach space of dimension $\geq 2$ over the real or
complex field ${\mathbb F}$ and
 ${\mathcal A}$  a standard operator algebra in ${\mathcal
B}(X)$. A map $\Phi:{\mathcal A} \rightarrow {\mathcal A}$ is said
to be  strong $3$-commutativity preserving if $[\Phi(A),\Phi(B)]_3 =
[A,B]_3$ for all $A, B\in{\mathcal A}$, where $[A,B]_3$ is the
3-commutator of $A,B$ defined by $[A,B]_3=[[[A,B],B],B]$. The main
result in this paper is shown that, if $\Phi$ is a surjective map on
${\mathcal A}$, then $\Phi$ is strong $3$-commutativity preserving
if and only if there exist
 a functional $h :{\mathcal A} \rightarrow {\mathbb F}$ and a scalar
$\lambda \in{\mathbb F}$ with $\lambda^4 = 1$ such that $\Phi(A) =
\lambda A + h(A)I$ for all $A \in{\mathcal A}$.

\end{abstract}
\maketitle

\section{Introduction}
Let ${\mathcal R}$ be a ring (or an algebra over a field ${\mathbb
F}$). Then ${\mathcal R}$ is a Lie ring (or Lie algebra) under the
Lie product $[A, B] = AB- BA$. Recall that a map $\Phi$ from
${\mathcal A}$ into itself is called a commutativity preserving map
if $[\Phi(A), \Phi(B)] = 0$ whenever $[A, B] = 0$ for all $A$,
$B\in{\mathcal R}$. The problem of characterizing commutativity
preserving maps had been studied
 intensively (see \cite{B, BC, MS, Z} and the references therein).
 Bell and Dail gave the conception of strong commutativity preserving maps in \cite{BD}. A
map $\Phi:{\mathcal R} \rightarrow {\mathcal R}$ is said to be
strong commutativity preserving if $[\Phi(A), \Phi(B)] = [A, B]$ for
any $A, B\in{\mathcal R}$. Clearly, strong commutativity preserving
maps must be commutativity preserving maps, but the inverse is not
true generally. The structure of linear (or additive) strong
commutativity preserving maps has been investigated in \cite{BD, BM,
DA, LL}. For nonlinear strong commutativity preserving maps, with
$\mathcal R$ being a prime unital ring containing a nontrivial
idempotent element, Qi and Hou in \cite{QXH} proved that every
nonlinear surjective strong commutativity preserving map
$\Phi:{\mathcal R} \rightarrow {\mathcal R}$  has the form $\Phi(A)
= \lambda A + f (A)$ for all $A \in{\mathcal R}$, where $\lambda
\in\{-1, 1\}$ and $f$ is a map from ${\mathcal R}$ into ${\mathcal
Z}_{\mathcal R}$, the center of ${\mathcal R}$. Particularly, this
result is true for maps on factor von Neumann algebras. In \cite{QH}
the nonlinear surjective strong commutativity preserving maps on
triangular algebras are studied. Recently, Liu in \cite{L} obtained
that a surjective strong commutativity preserving map
$\Phi:{\mathcal A} \rightarrow {\mathcal A}$ on von Neumann algebras
${\mathcal A}$ without central summands of type $I_1$ has the form
$\Phi(A) = ZA + f(A)$ for all $A \in{\mathcal A}$, where
$Z\in{\mathcal Z_{\mathcal A}}$ satisfies $Z^2 = I$ and $f$ is a map
from ${\mathcal A}$ into ${\mathcal Z_{\mathcal A}}$.

For a ring $\mathcal R$  and a positive integer $k$,  recall that
the $k$-commutator  of elements $A, B\in{\mathcal R}$ is defined by
$[A,B]_k = [[A,B]_{k-1},B]$ with $[A,B]_0 = A$ and $[A,B]_1 =[A,B]=
AB - BA$; a map $\Phi:{\mathcal R} \rightarrow {\mathcal R}$ is said
to be strong $k$-commutativity preserving if $
[\Phi(A),\Phi(B)]_k=[A,B]_k $ for all $A,B\in{\mathcal R}$.
Obviously, strong $k$-commutativity preserving maps are usual strong
commutativity preserving maps if $k=1$. It seems the study of the
problem of characterizing strong $k$-commutativity preserving maps
was started by
 \cite{Q}, in where it is shown that  a  nonlinear  surjective map on a unital prime ring containing a nontrivial idempotent
is strong $2$-commutativity preserving if and only if it has the
form $A\mapsto \lambda A+ f(A)$, where $\lambda$ is an element in
the extended centroid  of the ring satisfying $\lambda^3=1$ and $f$
is a map from the ring into its center. With $k$ increasing, the
problem of characterizing strong $k$-commutativity preserving maps
becomes much more difficult. Let $\mathcal R$ be a unital prime ring
containing a nontrivial idempotent. It is reasonable to conjecture
that a surjective map $\Phi:{\mathcal R}\to {\mathcal R}$ is strong
$k$-commutativity preserving if and only if there exist an element
$\lambda$ in the extended centroid of $\mathcal R$ with
$\lambda^{k+1}=1$ and a map $f:{\mathcal R}\to {\mathcal
Z}_{\mathcal R}$ such that $\Phi(A)=\lambda A+f(A)$ for all
$A\in{\mathcal R}$. However, we even do not know whether or not the
conjecture is true for $k=3$. The purpose of this paper is  to
answer the above conjecture affirmatively for the case when $k=3$
and the maps act on standard operator algebras on Banach spaces.

  Let ${\mathbb F} $ be the real field $\mathbb R$ or
the complex field $\mathbb C$ and $X$ be a Banach space over
$\mathbb F$. As usual, denote by ${\mathcal B}(X)$ the Banach
algebra of all bounded linear operators acting on $X$. Recall that,
a subalgebra ${\mathcal A}\subseteq{\mathcal B}(X)$ is called a
standard operator algebra if it contains the identity $I$ and all
finite rank operators.

The following is our main results.

\textbf {Theorem 1.1.} {\it Let $X$ be a  Banach space over the real
or complex field $\mathbb F$ with $\dim X\geq 2$ and ${\mathcal
A}\subseteq{\mathcal B}(X)$ be a standard operator algebra. Assume
that $\Phi:{\mathcal A} \rightarrow {\mathcal A}$ is a surjective
map. Then $\Phi$ is strong $3$-commutativity preserving if and only
if there exist a functional $h :{\mathcal A} \rightarrow {\mathbb
F}$ and a scalar $\lambda \in{\mathbb F}$ with $\lambda^4 = 1$ such
that $\Phi(A) = \lambda A + h(A)I$ for all $A \in{\mathcal A}$.}

\section{Proof of the main result}

 Before proving Theorem 1.1, we  give several  lemmas. In this
 section we always assume that $\mathcal A$ is a standard operator algebra on a real or complex
 Banach space $X$ and $\mathbb F$ stands for the real field $\mathbb
 R$ or the complex field $\mathbb C$, depending on $X$ is real or
 complex. For $x\in X$ and $f\in X^*$, we denote $x\otimes f$ for
 the rank one operator defined by $z\mapsto \langle z, f\rangle x$,
 where $\langle z, f\rangle=f(z)$. Note that, every operator of rank
 $\leq 1$ can be written in this form.

The first lemma is obvious by the main result in \cite{FS}.

 \textbf{Lemma 2.1.}  {\it Let ${\mathcal A}$ be a
 standard operator algebra, and let $A_i$, $B_i$, $C_j$, $D_j\in{\mathcal A}$, $i=1,2,\ldots, n$, such that $\sum{_{i=1}^n}
 A_{i}T
B_{i} = \sum{_{i=1}^n} C_{j}TD_{j}$ for all $ T\in {\mathcal A}$. If
$A_{1},\ldots,A_{n}$ are linearly independent, then each $B_{i}$ is
a linear combination of $D_{1},\ldots,D_{m}$. Similarly, if
$B_{1},\ldots,B_{m}$ are linearly independent, then each $A_{i}$ is
a linear combination of $C_{1},\ldots,C_{m}$. In particular, if
$ATB = BTA$ for all $T\in{\mathcal A}$, then $A$ and $B$ are linearly
dependent.}

\textbf{Lemma 2.2.} {\it Let ${\mathcal A}$ be a standard operator
algebra and $A \in{\mathcal A}$, if $[A, P]_3 = 0$ for any rank one
idempotent $P\in{\mathcal A}$, then $A \in \mathcal Z(\mathcal A) =
\{\lambda I:\lambda\in{\mathbb F}\}$.}

\textbf{Proof.} For any nonzero vector $x\in X$, there exists $f\in
X^*$ such that $\langle x, f \rangle = 1$. Let $P = x \otimes f$; then
$P$ is a rank one idempotent.  By the assumption,
$$0 = [A, P]_3 = [A, x \otimes f]_3 = [A, x \otimes f] = Ax \otimes
f - x \otimes A^*f,$$ which implies that $Ax $ is linearly dependent
of $x$. Thus, for any $x\in X$, there exists some scalar $\lambda_x
\in{\mathbb F}$ such that $Ax = \lambda_x x $. It follows that there
exists a scalar $\lambda \in{\mathbb F}$ such that $A = \lambda I$.
Hence $A \in \mathcal Z(\mathcal A) = \{\lambda I:\lambda\in{\mathbb
F}\}$. \hfill$\Box$

It was proved in \cite{Q} that, for $k = 2$,  if an element $S$
satisfies $[A,S]_k = 0$ for all $A$ in a prime ring, then $S$ is a
central element. But, for $k\geq 3$, this result is not true
anymore. For $k=3$, we have the following lemma, which reveals one
of the main difficulties to solve the problem of characterizing the
maps preserving strong $3$-commutativity.

 \textbf{Lemma 2.3.} {\it Let ${\mathcal A}$ be a
standard operator algebra, ${\mathcal N_2(\mathcal A)} = \{N \mid
N^2 = 0, N \in {\mathcal A}\}$ and $S \in{\mathcal A}$. Then
$[A,S]_3 = 0$ holds for all $A \in{\mathcal A}$ if and only if there
exist a scalar $\lambda \in{\mathbb F}$ and an element
$N\in{\mathcal N_2(\mathcal A)}$ such that $S = \lambda I + N$.}

{\bf Proof.} To check the ``if" part, assume that $S = \lambda I +
N$ with $N^2 = 0$. It is easily seen that
$$[A, S]_3 = [A, \lambda I + N]_3 = [A, N]_3 =AN^3-3NAN^2+3N^2AN-N^3A= 0$$
for any $ A\in{\mathcal A}.$

Next we check the `` only if " part.

Assume that $$[A, S]_3 = AS^3 - 3SAS^2 + 3S^2AS - S^3A = 0
\eqno(2.1)$$ for all $A \in{\mathcal A}$. By Lemma 2.1, we see that
$I$, $S$, $S^2$ and $S^3$ are linearly dependent.

If $I$ and $S$ are linearly dependent, then $S = \lambda I$ for some
$\lambda \in{\mathbb F}$ and the proof is done.

So we may suppose that $I$ and $S$ are linearly independent. We
assert that $I$, $S$ and $S^2$ are linearly dependent. In fact, if
$I$, $S$ and $S^2$ are linearly independent, then, by Lemma 2.1, we
get $S = \lambda I$, a contradiction. So there exist two scalars
$\alpha, \beta \in{\mathbb F}$ such that $S^2 = \alpha I + \beta S$.
Then
$$S^3 = S(\alpha I + \beta S) = \alpha S + \beta(\alpha I + \beta S)
= \alpha\beta I + (\alpha + \beta^2)S$$ and it follows from Eq.(2.1)
that
$$\begin{array}{rl}
0 = &AS^3 - 3SAS^2 + 3S^2AS - S^3A \\
= &A(\alpha\beta I + (\alpha + \beta^2)S) - 3SA(\alpha I + \beta S)
+ 3(\alpha I + \beta S)AS - (\alpha\beta I + (\alpha + \beta^2)S)A\\
= &\alpha\beta A + (\alpha + \beta^2)AS - 3\alpha SA - 3\beta SAS +
3\alpha I + 3\beta SAS - \alpha\beta A - (\alpha + \beta^2)SA
\\= &( 4\alpha + \beta^2)(AS - SA).
\end{array}$$
Since there  always exists $A \in{\mathcal A}$ such that $AS -
SA\neq 0$, we must have $4\alpha + \beta^2 = 0$.

Let $\lambda = \frac{\beta}{2}$, $N = S - \lambda I$, then
$\lambda^2 = - \alpha$, and
$$N^2 = S^2 - 2\lambda S + \lambda^2 I = S^2 - \beta S -\alpha I = 0.$$
So $S = \lambda I + N$ with $N \in {\mathcal N_2(\mathcal A)}$. The
proof of the Lemma 2.3 is completed. \hfill$\Box$

Now we are at the position to give our proof of the main theorem.

\textbf {Proof of Theorem 1.1.} The `` if " part of the theorem is
obvious. In the sequel, we always assume that
  $\Phi:{\mathcal A} \rightarrow {\mathcal A}$ is a surjective strong
$3$-commutativity preserving map. we will check the `` only if "
part by several steps.

\textbf {Step 1.} For any $A,B\in{\mathcal A}$, there exists a scalar
$\lambda_{A,B} \in{\mathbb F}$ such that $\Phi(A + B) = \Phi(A) +
\Phi(B) + \lambda_{A,B} I$.

 For any $A, B$ and $T\in\mathcal A$, we have
 \begin{align*}
 &[\Phi(A + B) - \Phi(A) - \Phi(B), \Phi(T)]_3\\
 &= [\Phi(A + B), \Phi(T)]_3 - [\Phi(A), \Phi(T)]_3 - [\Phi(B),
 \Phi(T)]_3\\
&= [A  + B, T]_3 - [A, T]_3 - [B, T]_3 = 0.
\end{align*}
By the surjectivity of $\Phi$ and Lemma 2.2, above equation implies
$\Phi(A + B) - \Phi(A) - \Phi(B)\in\{\lambda I:\lambda\in{\mathbb
F}\}$. So the assertion in Step 1 is true.

\textbf {Step 2.} $\Phi({\mathbb F}I) = {\mathbb F}I $ and
$\Phi({\mathbb F}I + {\mathcal N_2(\mathcal A)}) = {\mathbb F}I +
{\mathcal N}_2({\mathcal A})$.

 $\Phi({\mathbb F}I) = {\mathbb F}I $ is obvious by Lemma 2.2 and
 the surjectivity of $\Phi$.

If $N\in {\mathcal N}_2({\mathcal A})$ and $\lambda\in{\mathbb F}$,
then, for any $A\in{\mathcal A}$, we have
$$[\Phi(A), \Phi(\lambda I + N)]_3 = [A, \lambda I + N]_3
= [A, N]_3 = 0.$$ By the surjectivity of $\Phi$ and Lemma 2.3, we
see that $\Phi(\lambda I + N) \in {\mathbb F}I + {\mathcal
N}_2({\mathcal A})$.

On the other hand, if $\Phi(B) = \lambda I + N$ for some scalar
$\lambda$ and $N\in{\mathcal N}_2({\mathcal A})$, then
$$[A,B]_3 = [\Phi(A),\Phi(B)]_3 = [\Phi(A), \lambda I + N]_3 = [\Phi(A),
N]_3 = 0$$ for all $A\in{\mathcal A}$. By Lemma 2.3 again, we see
that $B \in {\mathbb F} I + {\mathcal N}_2({\mathcal A})$.

Hence $\Phi({\mathbb F} I + {\mathcal N_2(\mathcal A)}) = {\mathbb
F}I + {\mathcal N_2(\mathcal A)}$,
 completing the proof of the Step 2.

\textbf {Step 3.} For any nontrivial idempotent $P\in{\mathcal A}$,
there exist three scalars
 $\lambda_P$, $\mu_P, \mu_{I-P}\in{\mathbb F}$ and an element $N_P\in {\mathcal N_2(\mathcal
A)}$  such that
 $\Phi(P) = \lambda_P P + \mu_P I + N_P$ and $\Phi(I-P) = \lambda_P(I - P) + \mu_{I-P} I - N_{P}$.

Let $P\in{\mathcal A}$ be a nontrivial idempotent and  write $P_1 =
P$ and $P_2 = I - P$. Then $\mathcal A$ can be written as ${\mathcal
A} = {\mathcal A}_{11} + {\mathcal A}_{12} + {\mathcal A}_{21} +
{\mathcal A}_{22}$, where ${\mathcal A}_{ij} = P_i{\mathcal A}P_j\
  (i,j\in\{1, 2\})$. Write $\Phi(P) = S_{11} + S_{12}+ S_{21}+ S_{22}$ with
$S_{ij}\in{\mathcal A}_{ij}$.

We check Step 3 by several claims.

\textbf {Claim 3.1.} For any $A\in{\mathcal A}$, we have
$[A, \Phi(P)]_3\in{\mathcal A}_{12} + {\mathcal A}_{21}$.

Note that, for any $A\in{\mathcal A}$, we have
$$[A,P]_2 = P A (I - P) + (I - P) A P \in{\mathcal A}_{12} + {\mathcal
A}_{21}.$$ Also notice  that $[A, Q] = [A, Q]_3$ holds for any $A, Q\in
{\mathcal A}$ with $Q$ an idempotent. So we have $[A, P]_3 = [A, P]_5
= [[A, P]_3, P]_2 $, which implies that $$[\Phi(A), \Phi(P)]_3 =
[[\Phi(A), \Phi(P)]_3, P]_2$$ holds for any $A\in{\mathcal A}$.  By
the surjectivity of $\Phi$ one gets
$$[A, \Phi(P)]_3 = [[A, \Phi(P)]_3, P]_2\ \ \mbox{\rm for all}\ A\in{\mathcal A}.$$
Therefore, for any $A$ we have $$[A, \Phi(P)]_3\in{\mathcal A}_{12} +
{\mathcal A}_{21}.$$

The next claim is one of the key steps for our proof.

\textbf{Claim 3.2.} $S_{12} = S_{21} = 0.$

Taking $A = A_{11}\in{\mathcal A}_{11}$ and applying Claim 3.1, we
have
$$\begin{array}{rl}
 0=&(I - P) [A_{11},\Phi(P)]_3 (I - P)\\
  = &3S_{21}S_{11}A_{11}S_{12}+3S_{22}S_{21}A_{11}S_{12} - 3
S_{21}A_{11}S_{11}S_{12} - 3S_{21}A_{11}S_{12}S_{22},
\end{array}$$
 which implies that
  $$(S_{21}S_{11} +
S_{22}S_{21})A_{11}S_{12} = S_{21}A_{11}(S_{11}S_{12} +
S_{12}S_{22})$$ holds for all $A_{11}\in{\mathcal A}_{11}$,
 that is,
 $$(S_{21}S_{11} + S_{22}S_{21})P A P S_{12}
= S_{21}P A P(S_{11}S_{12} + S_{12}S_{22})\eqno(2.2)$$ holds for all
$A\in{\mathcal A}$.
 Thus, by Lemma 2.1,
$S_{21}S_{11} + S_{22}S_{21}$ and $S_{21}$ must be linearly
dependent. It follows that
$$S_{21}S_{11} + S_{22}S_{21} = \lambda_1S_{21}\eqno(2.3)$$ for some
$\lambda_1\in{\mathbb F}$. Again by Eq.(2.2), we get
 $$S_{11}S_{12} + S_{12}S_{22} = \lambda_1 S_{12}.\eqno(2.4)$$
 Taking $ A = A_{12}\in{\mathcal A}_{12}$ and applying Claim 3.1, one
gets
\begin{align*}
  0 = &P [A_{12}, \Phi(P)]_3 P \\
 = &A_{12}S_{21}S_{11}^2 + A_{12}S_{21}S_{12}S_{21}
 + A_{12}S_{22}S_{21}S_{11} + A_{12}S_{22}S_{21}S_{12}\\
  &- 3S_{11}A_{12}S_{21}S_{11} - 3S_{11}A_{12}S_{22}S_{21} +
3S_{11}^2A_{12}S_{21} + 3S_{12}S_{21}A_{12}S_{21}
\end{align*}
and
\begin{align*}
0 = & (I - P) [A_{12},\Phi(P)]_3 (I - P) \\
=& -S_{21}S_{11}^2A_{12} - S_{22}S{21}S_{11}A_{12}
 - S_{21}S_{12}S_{21}A_{12} - S_{22}^2S_{21}A_{12}\\
&- 3S_{21}A_{12}S_{21}S_{12} - 3S_{21}A_{12}S_{22}^2 +
3S_{21}S_{1}1A_{12}S_{22} + 3S_{22}S_{ 21}A_{12}S_{22}.
\end{align*}
Thus
$$\begin{array}{rl}
& A_{12}[S_{21}(S_{11}^2 + S_{12}S_{21})
 + S_{22}(S_{21}S_{11} + S_{22}S_{21})] \\
 = & 3S_{11}A_{12}(S_{21}S_{11} + S_{22}S_{21})
 - 3(S_{11}^2 + S_{12}S_{21})A_{12}S_{21}
 \end{array} \eqno(2.5)$$
  and
 $$\begin{array}{rl} &[S_{21}(S_{11}^2 + S_{12}S_{21})
 + S_{22}(S_{21}S_{11} + S_{22}S_{21})]A_{12}\\
 = & 3S_{21}A_{12}(S_{22}^2 + S_{21}S_{12})
 - 3(S_{21}S_{11} + S_{22}S_{21})A_{12}S_{22}
 \end{array} \eqno(2.6)$$
hold for all $A_{12}\in{\mathcal A}_{12}$.

Similarly, by taking $ A = A_{21}$ and applying Claim 3.1, one gets
\begin{align*}0
 =& P [A_{21},\Phi(P_1)]_3 P \\
 =& S_{11}^2S_{12}A_{21} + S_{11}S_{12}S_{22}A_{21}
 + S_{12}S{22}^2A_{21} + S_{12}S_{21}S_{12}A_{21}\\
&- 3S_{11}S_{12}A_{21}S_{11} - 3S_{12}S_{22}A_{21}S_{11}
 + 3S_{12}A_{21}S_{11}^2 + 3S_{12}A_{21}S_{12}S_{21}
\end{align*}
and
\begin{align*}
0 =& (I - P)[A_{21},\Phi(P)]_3 (I - P) \\
=& A_{21}S_{11}^2S_{12} + A_{21}S_{11}S_{12}S_{22}
 + A_{21}S_{12}S_{22}^2 +  A_{21}S_{12}S_{21}S_{12} \\
&- 3S_{22}A_{21}S_{11}S_{12} - 3S_{22}A_{21}S_{12}S_{22} +
3S_{22}^2A_{21}S_{12} + 3S_{21}S_{12})A_{21}S_{12},
\end{align*}
which imply that
$$\begin{array}{rl} &[S_{11}(S_{11}S_{12} + S_{12}S_{22})
 + S_{12}(S_{22}^2 + S_{21}S_{12})]A_{21} \\ = &3(S_{11}S_{12} + S_{12}S_{22})A_{21}S_{11}
 - 3S_{12}A_{21}(S_{11}^2 + S_{12}S_{21})\end{array}\eqno(2.7)$$
and
$$\begin{array}{rl} & A_{21}[S_{11}(S_{11}S_{12} + S_{12}S_{22})
 + S_{12}(S_{22}^2 + S_{21}S_{12})] \\= & 3S_{22}A_{21}(S_{11}S_{12} + S_{12}S_{22})
 - 3(S_{22}^2 + S_{21}S_{12})A_{21}S_{12}\end{array} \eqno(2.8)$$ hold for all $A_{21}\in{\mathcal A}_{21}$.
Obviously, Eqs.$(2.3)\thicksim(2.8)$ together are equivalent to say
that
$$\begin{cases}
PA[S_{21}(S_{11}^2 + S_{12}S_{21}) + S_{22}(S_{21}S_{11} +
S_{22}S_{21})]
 = [3\lambda_1 S_{11} - 3(S_{11}^2 + S_{12}S_{21})]AS_{21},\\
[S_{21}(S_{11}^2 + S_{12}S_{21}) + S_{22}(S_{21}S_{11} +
S_{22}S_{21})]A (I - P)
= S_{21}A[3\lambda_1 S_{22}-3(S_{22}^2 + S_{21}S_{12})],\\
[S_{11}(S_{11}S_{12} + S_{12}S_{22}) + S_{12}(S_{22}^2 +
S_{21}S_{12})]A P
 = S_{12}A P [3\lambda_1S_{11} - 3(S_{11}^2 + S_{12}S_{21})],\\
(I - P) A [S_{11}(S_{11}S_{12} + S_{12}S_{22}) + S_{12}(S_{22}^2 +
S_{21}S_{12})]\\
 \quad =[3\lambda_1S_{22}
 - 3(S_{22}^2 + S_{21}S_{12})](I - P)A S_{12}
\end{cases}$$
hold for all $A\in{\mathcal A}$.

Thus, by Lemma 2.1, we see that $3 \lambda_1S_{11} - 3(S_{11}^2 +
S_{12}S_{21})$ and $P$ are linearly dependent. So
$$3\lambda_1S_{11} - 3(S_{11}^2 + S_{12}S_{21}) = \lambda_2 P \eqno(2.9)$$
for some $\lambda_2\in{\mathbb F}$. This entails that
$$[S_{21}(S_{11}^2 + S_{12}S_{21}) + S_{22}(S_{21}S_{11} + S_{22}S_{21})]
 = \lambda_2S_{21},\eqno(2.10)$$
 $$[S_{11}(S_{11}S_{12} + S_{12}S_{22}) + S_{12}(S_{22}^2 + S_{21}S_{12})]
 = \lambda_2S_{12} \eqno(2.11)$$
 and
 $$[3\lambda_1S_{22} - 3(S_{22}^2 + S_{21}S_{12})]
 = \lambda_2(I - P). \eqno(2.12)$$
 It is easily checked by Eqs.$(2.9)\thicksim(2.12)$ that
 $$(\lambda_1^2 - \frac{4\lambda_2}{3})S_{21} = 0\quad{\rm and}\quad
  (\lambda_1^2 - \frac{4\lambda_2}{3})S_{12} = 0.$$

 If $\lambda_1^2 - \frac{4\lambda_2}{3} = 0$, then we have
 $$\begin{array}{rl}&(\Phi(P)-\frac{\lambda_1}{2}I)^2 \\
 =& (S_{11} + S_{12} + S_{21} + S_{22} - \frac{\lambda_1}{2}I)^2\\
 =&(S_{11}^2 + S_{12}S_{21}) + (S_{11}S_{12} + S_{12}S_{22}) + (S_{21}S_{11} + S_{22}S_{21}) \\
 &+ (S_{22}^2 + S_{21}S_{12}) - \lambda_1S_{11} - \lambda_1S_{12} - \lambda_1S_{21} - \lambda_1S_{22} + \frac{\lambda_1^2}{4}I\\
 =&(\lambda_1S_{11} - \frac{\lambda_2}{3}P) + \lambda_1S_{12} + \lambda_1S_{21}
 + (\lambda_1S_{22} - \frac{\lambda_2}{3}(I - P)) - \lambda_1S_{11} \\ & - \lambda_1S_{12} - \lambda_1S_{21} - \lambda_1S_{22} + \frac{\lambda_1^2}{4}I\\
 =&(\frac{\lambda_1^2}{4} - \frac{\lambda_2}{3})I = 0.
 \end{array}.$$
 Let $N_P = \Phi(P) - \frac{\lambda_1}{2}I$; then $\Phi(P) = \frac{\lambda_1}{2}I + N_P$ with $N_P^2 = 0$. But, by
 Step 2, this implies  $P = \lambda I + N$ for some $\lambda \in{\mathbb F}$ and some $
 N \in {\mathcal N_2(\mathcal A)}$, a contradiction. So we must have $\lambda_1^2 -
\frac{4\lambda_2}{3} \neq 0$. Hence $S_{12} = S_{21} = 0$ and
$\Phi(P) = S_{11} + S_{22}$.

  Next, let us determine $S_{11}$ and $S_{22}$.

 Taking $ A = A_{11}$ and applying Claim 3.1, we have
 $$[A_{11}, S_{11}]_3 =A_{11}S_{11}^3 - 3 S_{11}A_{11}S_{11}^2
 + 3S_{11}^2A_{11}S_{11} - S_{11}^3A_{11} =P[A_{11}, \Phi(P)]_3P =   0$$
holds for all $A_{11}\in{\mathcal A}_{11}$. As ${\mathcal A}_{11}$
is clearly a standard operator algebra in ${\mathcal B}(PX)$, by
Lemma 2.3, $S_{11} = \lambda_1P +
 N_1$ with $N_1^2 = 0$ for some $\lambda_1 \in{\mathbb F}$ and some $N_1\in{\mathcal N_2(\mathcal A_{11})}$.
  A similar argument can show  that $S_{22} = \mu_P (I - P) +
 N_2$ for some $\mu_P\in{\mathbb F}$ and some $N_2\in{\mathcal N_2(\mathcal A_{22})}$. Hence
  $$\Phi(P) = \lambda_1P + N_1 + \mu_P (I - P) +
 N_2 = \lambda_P P + \mu_P I + N_P,$$
 where $\lambda_P = \lambda_1 - \mu_P$ with $N_P = N_1 + N_2$, $N_P^2 =
 0$ and $N_1N_2=N_2N_1=0$.

 Similarly, we can get $\Phi(I - P) = \lambda_{I - P}(I - P) + \mu_{I - P} I + N_{I -P}$ for some $\lambda_{I-P},\mu_{I-P}\in{\mathbb F}$
 and $N_{I-P}\in{\mathcal N}_2({\mathcal A})$.

\textbf{Claim 3.3.} $\lambda_P = \lambda_{(I - P)}$ and $N_{I - P} =
- N_P$.

By Step 1 and Claim 3.2, there exists some scalar $c\in{\mathbb F}$, such
that
$$\begin{array}{rl}
\Phi(I) &= \Phi(P) + \Phi(I - P) + cI\\
 &= \lambda_P P + \mu_P I + N_P + \lambda_{I - P}(I - P) + \mu_{I - P} I +
 N_{I - P} + cI\\
 &= (\lambda_P - \lambda_{I - P})P + (\mu_P¡¡+ \mu_{I - P} + c) I +  N_P + N_{I -
 P}.
 \end{array}\eqno(2.13)$$
Since $\Phi(I)\in{\mathbb F} I$ by Step 2,  there exist two scalars
$\alpha, t \in{\mathbb F}$, with $t = \lambda_P - \lambda_{I - P}$
such that
$$tP + N_P + N_{I - P} = \alpha I. \eqno(2.14)$$
 So $$A(tP + N_P + N_{I - P})
=(tP + N_P + N_{I - P})A \ \ \mbox{\rm for\ all}\  A\in{\mathcal A}.$$
 In particular, taking $A = N_P$ gives
 $$tN_P P + N_PN_{I - P} = tP N_P + N_{I - P}N_P. \eqno(2.15)$$
Recall that $N_P = N_1 + N_2$ with $N_i\in{\mathcal N}_2({\mathcal
A}_{ii})$, $i = 1,2$. So $tN_P P = tP N_P = tN_1$ and hence,  by
Eq.(2.15), we have $N_P N_{I - P} = N_{I - P}N_P$. Then, it follows from
Eq.(2.14) that
$$(tP - \alpha I)^2 = (N_P + N_{I - P})^2 = 2N_P N_{I - P}. \eqno(2.16)$$
Note that $(N_P N_{I - P})^2 = 0$ as $N_P N_{I - P} = N_{I - P}N_P$. Thus,
Eq.(2.16) entails that $t = \alpha = 0$, which, by Eq.(2.14), implies
that $\lambda_{I - P} = \lambda_P$ and $N_{I - P} = -N_P$, as desired.

\if that is
$$(t^2 - 2\alpha t) P = 2N_p N_{I-p} - \alpha^2 I\eqno(2.15)$$
Assume $t^2 - 2\alpha t \neq 0$,then
 $$P = \frac{2}{t^2 - 2\alpha t} N_PN_{I-P} - \frac{\alpha^2}{t^2 - 2\alpha t} I \eqno(2.16)$$
since $(N_P N_{I-P})^2 = N_P N_{I-P}^2 N_P = 0$, so Eq.(2.16) is a
contradiction. Thus $t^2 - 2\alpha t = t(t - 2\alpha) = 0$, this is
$ t = 0$ or $t - 2\alpha = 0$. Assume $ t \neq 0$, thus $t - 2\alpha
= 0$. By Eq.(2.14), we have
$$\alpha P - \alpha I = N_P + N_{I - P}.$$
Taking twice square from both side above equation, we get $\alpha^4
I = 0,$ a contradiction, then $t = 0$, So $\lambda_P = \lambda_{(I -
P)}$.

Putting $t = 0$ into Eq.(2.14), one gets
 $$N_P + N_{I - P} = \alpha I$$
Taking twice square from both side above equation, we get $\alpha^4
I = 0$, which implies $\alpha = 0$, that is $N_{I - P} = - N_P$.\fi

\textbf {Step 4.}  For any  idempotent $P\in{\mathcal A}$  of rank
$\leq 2$, we have $N_P = 0$; For any rank-1 idempotents $P, Q$ with
$PQ = QP = 0$, we have  $\lambda_P = \lambda_Q$.

 Assume that $P$ is a rank-1 idempotent.  By Step 3, there
exist two scalars $\lambda_P, \mu_P\in{\mathbb F}$ and an element
$N_P \in {\mathcal N_2(\mathcal A)}$, such that $\Phi(P)= \lambda_P
P + \mu_P I + N_P$. Taking the space decomposition $X = X_1\dot{+}
X_2$ so that $P$ has the matrix representation
 $P = \left(
\begin{array}{cc}
1&0\\
0&0\\
\end{array}
\right)$. Since $N_P = N_1 + N_2$ with $N_1\in P{\mathcal A}P$ and
$N_2\in(I - P){\mathcal A}(I - P)$,
 we must have $N_1 = 0$ and hence
 $N_P = \left(
  \begin{array}{cc}
  0&0\\
0&S\\
\end{array}
\right)$
 with $S^2 = 0$.

Assume $S \neq 0$, then there exists  rank-1 idempotent
${P_2^\prime} \in B((I - P)X)$, such that ${P_2^\prime}S \neq 0$. Let
$P_2 = 0 \oplus{P_2^\prime}$, then $P P_2 = P_2 P$ and with respect to the
space decomposition $X = PX \dot{+} P_2X \dot{+} (I - P)(I - P_2)X$,
we have
$$P = \left( \begin{array}{cccc}1&0&0\\
0&0&0\\0&0&0\\ \end{array} \right) \quad {\rm and}\quad P_2 = \left( \begin{array}{cccc}0&0&0\\
0&1&0\\0&0&0\\ \end{array} \right).$$
Let $$Q =P+P_2= \left( \begin{array}{cccc}1&0&0\\
0&1&0\\0&0&0\\ \end{array} \right).$$ Obviously, $Q^2 = Q$. By Step
1 and Step 3, there exist scalars $c,\lambda_P, \lambda_{P_2}$,
$\lambda_Q$, $\mu_Q$, $\mu_{P}$, $\mu_{P_2}\in{\mathbb F}$ and
elements $N_Q, N_P, N_{P_2}\in {\mathcal N_2(\mathcal A)}$
 with
 $$N_Q = \left(\begin{array}{ccc}t_{11}&t_{12}&0\\
 t_{21}&t_{22}&0\\0&0&T\\\end{array}\right),\
N_P =
\left(\begin{array}{ccc}0&0&0\\0&s_{22}&S_{23}\\0&S_{32}&S_{33}\\\end{array}
\right)\ {\rm and}\ N_{P_2} =
\left(\begin{array}{ccc}r_{11}&0&R_{13}\\0&0&0\\R_{31}&0&R_{33}\\\end{array}\right)$$
 such that
$$\begin{array}{rl} \lambda_Q Q + \mu_Q I + N_Q = & \Phi(Q) = \Phi(P_1) +
\Phi(P_2) + cI \\ = &\lambda_P P + \mu_P I + N_P + \lambda_{P_2} P_2
+ \mu_{P_2} I + N_{P_2}. \end{array}$$  So we have
$$\begin{array}{rl} & \left(
\begin{array}{ccc}
\lambda_Q + \mu_Q&0&0\\
0&\lambda_Q + \mu_Q&0\\
0&0&\mu_Q\\
 \end{array}
 \right)
+ \left(
 \begin{array}{ccc}
 t_{11}&t_{12}&0\\
t_{21}&t_{22}&0\\
0&0&T\\
\end{array}
\right)\\
 = & \left(
\begin{array}{ccc}
\lambda_P + \mu_P + \mu_{P_2} + c&0&0\\
0&\lambda_{P_2} + \mu_P + \mu_{P_2} + c&0\\
0&0&\mu_P + \mu_{P_2} + c\\
\end{array}
\right)
+ \left(
 \begin{array}{ccc}
 r_{11}&0&R_{13}\\
0&s_{22}&S_{23}\\
R_{31}&S_{32}&S_{33}+R_{33}\\
\end{array}
 \right), \end{array}$$
which entails that $t_{12} = t_{21} = 0$, $R_{13} = 0$, $R_{31} = 0$,
$S_{23} = 0$ and $ S_{32} = 0$. Since $N_Q^2 = 0$, we get also
$t_{11} = t_{22} = 0$. Then, as $N_P^2 = 0$ and $N_{P_2}^2 = 0$, one
gets further that  $s_{22} =  r_{11} = 0$. So $N_P$ has the form
 $N_P = \left(
 \begin{array}{ccc}
 0&0&0\\
0&0&0\\
0&0&S_{33}\\
 \end{array}
 \right)$;
that is $S =\left(
\begin{array}{cc}
0&0\\
0&S_{33}\\
 \end{array}
 \right)$. However, as
 ${P_2}' =
  \left(
  \begin{array}{cc}
   1&0\\
0&0\\
\end{array}
\right)$, we get a contradiction ${P_2}'S = S{P_2}' = 0$. Hence we
must have
 $N_P = 0$. Moreover,
 $$\lambda_P + \mu_P + \mu_{P_2} + c = \lambda_Q + \mu_Q = \lambda_{P_2} + \mu_P + \mu_{P_2} +
 c.$$ This forces that $\lambda_{P_2} = \lambda_P$. It follows that $T = 0$;
 that is $N_Q = 0$, too.

Note that, for the case when $\dim X\geq 3$, we also have
$\mu_Q = \mu_P + \mu_{P_2} + c$ and hence $\lambda_Q = \lambda_P$.

 If $P$ is an idempotent of rank 2, taking rank-1 idempotent
 $P_1, P_2$ so that $P = P_1 + P_2$, then the above argument shows that
 $N_P = 0$. So, the assertion of Step 4 is true.

\if By now we have shown that, for any  idempotent $P$ of rank $\leq
2$ we have
 $N_P = 0$, and for any rank-1 idempotents $P, Q$ with $PQ = QP = 0$ we
 have $\lambda_P = \lambda_Q$.\fi

\textbf {Step 5.} Assume $\dim X\geq 3$. $\lambda^4 = 1$ and there
exists a functional $h :{\mathcal A} \rightarrow {\mathbb F} $ such
that $\Phi(A) = \lambda A + h(\mathcal A)I$ for any $A\in{\mathcal
A}$.

We first show that, if $\dim X\geq 3$, then there exists a scalar
$\lambda\in {\mathbb F}$ such that, for any rank-1 idempotent
$P\in{\mathcal A}$, $\Phi(P) = \lambda P + \mu_P I$ for some
$\mu_P\in{\mathbb F}$. In fact, for any rank one idempotent
operators $P, Q$, as $\dim X\geq 3$,
 there exists rank-1 idempotent $R$ such that $PR = RP = QR = RQ =0$. It
follows from Step 4 that $\lambda_P = \lambda_R = \lambda_Q$. Hence,
there exists a scalar $\lambda$ such that $\lambda_P = \lambda$ for
any rank-1 idempotent $P$.

Thus, for any rank-1 idempotents $P$, $Q \in{\mathcal A}$,   we have
$$[P, Q]_3 = [\Phi(P), \Phi(Q)]_3 = [\lambda P, \lambda Q]_3 = \lambda^4 [P,
Q]_3.$$
 It follows that
$\lambda^4 = 1$ as there are rank-1 idempotents $P,Q$ so that
$[P, Q]_3\not = 0$.

Now,  for any $A$ and any rank-1 idempotent $P$ in ${\mathcal A}$,
we have
$$[A, P]_3 = [\Phi(A), \Phi(P)]_3 = [\Phi(A), \lambda P]_3 = \lambda^3 [\Phi(A),
 P]_3,$$
 which is equivalent to $[\Phi(A) - \lambda A, P]_3 = 0$
holds for any rank-1 idempotent $P\in{\mathcal A}$.   By applying
Lemma 2.2, there exists a scalar $\mu_A \in{\mathbb F}$, such that
$\Phi(A) - \lambda A = \mu_A I$, that is
$$\Phi(A) = \lambda A + \mu_A I.$$

Define a functional $h :{\mathcal A} \rightarrow {\mathbb F}$  by
$h(A) = \mu_A$ for any  $A \in{\mathcal A}$.
 Then we get $\Phi(A) = \lambda A + h(A)I$. This completes the proof
 Theorem 1.1 for the case when $\dim X\geq 3$.

\textbf {Step 6.}  Assume $\dim X = 2$. There exists a scalar
$\lambda$ with $\lambda^4 = 1$ and a functional $h:{\mathcal A}\to
{\mathbb F}$ such that $\Phi(A) = \lambda A + h( A)I$ for any
$A\in{\mathcal A}$.

Fixing a basis  $\{e_1, e_2\}$ of $X$, ${\mathcal A}$ may be
identified with the matrix algebra $ {\mathcal M}_2={\mathcal
M}_2(\mathbb F)$ of all 2 by 2 matrices over $\mathbb F$. So, we may
regard the map $\Phi$ as a map  between ${\mathcal M}_2(\mathbb F)$
in the case when $\dim X = 2$.

Let $E_{11} = \left(
\begin{array}{cc}
  1 &0 \\
  0&0
\end{array}
\right)$, $ E_{22} = \left(
\begin{array}{cc}
  0&0 \\
  0&1
\end{array}
\right)$, $ E_{12} = \left(
\begin{array}{cc}
  0&1 \\
  0&0
\end{array}
\right)$ and $ E_{21} = \left(
\begin{array}{cc}
  0&0 \\
 1&0
\end{array}
\right)$. Then every $A\in {\mathcal A}$ can be written as $  A =
a_{11}E_{11} + {a}_{12}E_{12} + {a}_{21}E_{21} + {a}_{22}E_{22}$,
where ${a}_{ij} \in{\mathbb F}$.

By Step 3-4, there exist  scalars $\lambda, \mu_i$ such that
$\Phi(E_{ii}) = \lambda E_{ii} + \mu_i I$, $i = 1, 2$.

 We will check Step 6 by three claims.

 \textbf {Claim 6.1.} If $1\leq i\not = j\leq 2$, then, for any $a_{ij}\in\mathbb F$, there
exists a scalar $\mu_{a_{ij}}\in{\mathbb F}$, such that
$\Phi(a_{ij}E_{ij}) = \lambda^{-3}a_{ij}E_{ij} + \mu_{a_{ij}} I.$

Here, we only give the proof for the case when $(i,j) = (1,2)$. The
proof for $(i,j) = (2,1)$ is similar.

For any $a_{12}\in{\mathbb F}$, write $\Phi(a_{12}E_{12}) =
s_{11}E_{11} + s_{12}E_{12} + s_{21}E_{21} + s_{22}E_{22}$ with
$s_{ij}\in{\mathbb F}$. Then, by Step 4, we see that
$$\begin{array}{rl}
 - a_{12}E_{12} &= [a_{12}E_{12}, E_{11}]_3 = [\Phi(a_{12}E_{12}), \Phi(E_{11})]_3\\
 &= [s_{11}E_{11} + s_{12}E_{12} + s_{21}E_{21}+ s_{22}E_{22}, \lambda E_{11} + \mu_1 I ]_3 \\
 &= \lambda^3 [s_{11}E_{11} + s_{12}E_{12} + s_{21}E_{21} + s_{22}E_{22}, E_{11}]_3 \\
 &= \lambda^3(s_{21}E_{21} - s_{12}E_{12}),
\end{array}$$
which entails that $s_{21} = 0$ and $\lambda^3 s_{12} = a_{12}$. So
$s_{12} = \lambda^{-3} a_{12}$ and
$$\Phi(a_{12}E_{12}) = \lambda^{-3} a_{12}E_{12} + s_{11}E_{11} + s_{22}E_{22}.$$
On the other hand, by Step 2, we get $\Phi(a_{12}E_{12}) = \mu I +
N$ for some $\mu\in{\mathbb F}$ and $N\in{\mathcal N_2({\mathcal
M}_2)}$ with $N = \left(
\begin{array}{cc}
  r_{11} &r_{12}\\
  r_{21}&r_{22}
\end{array}
\right)$. It follows that
$$\begin{array}{rl}
&\lambda^{-3} a_{12}E_{12} + s_{11}E_{11} + s_{22}E_{22} = \mu I + N \\
&= (\mu + r_{11})E_{11} + r_{12}E_{12} + r_{21}E_{21} + (\mu +
r_{22})E_{22}.
\end{array}$$
\if that is
 $$\left(
\begin{array}{cc}
  s_{11}&\lambda^{-3} a_{12}\\
  0&s_{22}
\end{array}
\right) = \left(
\begin{array}{cc}
 \mu + r_{11}&r_{12}\\
  r_{21}&\mu + r_{22}
\end{array}
\right),$$\fi This gives that $r_{21} = 0$, $r_{12} = \lambda^{-3}
a_{12}$. Since $N^2 = 0$,
 we must have
$r_{11} = r_{22} = 0$, and then $s_{11} = s_{22} = \mu$.
 Thus one gets
 $$\Phi(a_{12}E_{12}) = \lambda^{-3} a_{12}E_{12} + \mu_{a_{12}} E_{11} + \mu_{a_{12}} E_{22} = \lambda^{-3} a_{12}E_{12} + \mu_{a_{12}}I,$$
 with $\mu_{a_{12}} = \mu$, as desired.

\textbf {Claim 6.2.} $\lambda^4 = 1$ and, for any $a_{ii}\in{\mathbb
F}_{ii}$, $(i\in\{1,2\})$, there exists a scalar
$\mu_{a_{ii}}\in{\mathbb F}$, such that $\Phi(a_{ii}E_{ii}) =
\lambda a_{ii}E_{ii} + \mu_{a_{ii}} I.$

Still, we only prove that the claim holds for the case $i=1$.

Take any nonzero $a_{11}\in{\mathbb F}$ and write
$\Phi(a_{11}E_{11}) = s_{11}E_{11} + s_{12}E_{12} + s_{21}E_{21}+
s_{22}E_{22}$ with $s_{ij}\in{\mathbb F}$. By Step 4, we obtain that
$$\begin{array}{rl}
 0 &= [a_{11}E_{11}, E_{11}]_3 = [\Phi(a_{11}E_{11}), \Phi(E_{11})]_3\\
 &= [s_{11}E_{11} + s_{12}E_{12} + s_{21}E_{21}+ s_{22}E_{22}, \lambda E_{11} + \mu_1 I ]_3 \\
 & = \lambda^3(s_{21}E_{21} - s_{12}E_{12}).
\end{array}$$
Hence $s_{21} = s_{12} = 0$ and $\Phi(a_{11}E_{11}) = s_{11}E_{11} +
s_{22}E_{22}$. By Claim 6.1, there exists a scalar
$\mu_{E_{12}}\in{\mathbb F}$, such that $\Phi(E_{12})=\lambda^{-3}
E_{12}+\mu_{E_{12}}I$. Thus we have
$$\begin{array}{rl}
 - a_{11}^3E_{12} &= [E_{12}, a_{11}E_{11}]_3 = [\Phi(E_{12}), \Phi(a_{11}E_{11})]_3\\
 &= [\lambda^{-3} E_{12} + \mu_{E_{12}} I, s_{11}E_{11} + s_{22}E_{22} ]_3 \\
 & = \lambda^{-3}(s_{22} - s_{11})^3E_{12}.
\end{array}$$
So
$$(s_{11} - s_{22})^3 = \lambda^3 a_{11}^3,$$
which implies that
$$s_{11} - s_{22} = \delta \lambda a_{11} $$ for some $\delta$ with
$\delta^3=1$.
 It follows that
$$\Phi(a_{11}E_{11}) = s_{11}E_{11} + s_{22}E_{22} = (\delta \lambda a_{11} + s_{22})E_{11} + s_{22}E_{22}= \delta \lambda a_{11}E_{11} + s_{22} I.$$
Applying Step 1 and Claim 6.1, it is easily checked that
$$\begin{array}{rl} 4a_{11}(E_{12} - E_{21})
&= [a_{11}E_{11}, E_{12} + E_{21}]_3 \\
&= [\Phi(a_{11}E_{11}), \Phi(E_{12} + E_{21})]_3
 = [\delta \lambda a_{11}E_{11}, \lambda^{-3}(E_{12} + E_{21})]_3 \\
 &= \delta \lambda^{-8} a_{11}[E_{11}, E_{12} + E_{21}]_3
 = 4\delta \lambda^{-8} a_{11}(E_{12} - E_{21}).
\end{array}$$
So we get $\delta = \lambda^8$ and then
$$\Phi(a_{11}E_{11}) =  \lambda^9 a_{11}E_{11} + \mu_{a_{11}} I $$
with $\mu_{a_{11}} = s_{22}$. Note that, we already had
$\Phi(E_{11}) = \lambda E_{11} + \mu_1 I$. This gives $\lambda^8=1$ and
hence
$$\Phi(a_{11}E_{11}) = \lambda a_{11}E_{11} + \mu_{a_{11}} I $$
holds for all $a_{11}\in{\mathbb F}$. In fact, we have  $\lambda^4 =
1$. To see this, notice that,
  by Step 1 and Claim 6.1, we have
$$\begin{array}{rl}
 E_{21} + E_{22} - E_{11} &= [E_{21}, E_{11} + E_{12}]_3 = [\Phi(E_{21}), \Phi(E_{11} + E_{12})]_3 \\
&= [\lambda^{-3}E_{21}, \lambda E_{11} + \lambda^{-3}E_{12}]_3 =
E_{21} + \lambda^{-4}E_{22} - \lambda^{-4}E_{11},
\end{array}$$
which implies that $\lambda^{4} = 1$.

\textbf {Claim 6.3.} There exists a functional $h$ on ${\mathcal
M}_2$ such that $\Phi(A) = \lambda A + h(\mathcal A)I$ for any
$A\in{\mathcal M}_2$.

 As $\lambda ^{-3}=\lambda$, by Claims 6.1-6.2, we have  that $\Phi(a_{ij}E_{ij})= \lambda a_{ij}E_{ij} + \mu_{a_{ij}} I$
  for any $i, j\in\{1, 2\}$ and any $a_{ij}\in {\mathbb F}$.
Then, for any $A\in{\mathcal M}_2$, writing $A = a_{11}E_{11} +
a_{12}E_{12} + a_{21}E_{21} + a_{22}E_{22}$ with $a_{ij}\in {\mathbb
F}$ and applying  Step 1,  there exists some scalar $c_A$ such that
$$\begin{array}{rl}
 \Phi(A) = & \Phi(a_{11}E_{11}) + \Phi(a_{12}E_{12}) + \Phi(a_{21}E_{21})+ \Phi(a_{22}E_{22}) + c_A I \\
=& \lambda a_{11}E_{11} + \mu_{a_{11}} I + \lambda a_{12}E_{12} +
\mu_{a_{12}} I + \lambda a_{21}E_{21} \\ &+ \mu_{a_{21}} I + \lambda
a_{22}E_{22} + \mu_{a_{22}} I + c_AI\\ = & \lambda A + (\mu_{a_{11}}
+ \mu_{a_{12}} + \mu_{a_{21}} + \mu_{a_{22}} + c_A) I = \lambda A +
\mu_A I.
\end{array}$$
Let $h :{\mathcal A} \rightarrow {\mathbb F}$ be a functional
defined  by $h(A) = \mu_A$ for any  $A \in{\mathcal M}_2$.
 Then we get $\Phi(A) = \lambda A + h(A)I$ for any $A\in{\mathcal
 M}_2$.

 The proof of Theorem 1.1 is completed. \hfill$\Box$

%{\bf Acknowledgement.} The authors wish to give their thanks to the
%referees for helpful comments and suggestions.


\begin{thebibliography}{99}


\bibitem{BD} H. E. Bell, M. N. Daif, On commutativity and strong commutativity preserving maps, Can. Math. Bull., 37 (1994), 443-447.


\bibitem{B} M. Bre\v{s}ar. Commuting traces of biadditive mappings, commutativity preserving mappings and Lie mappings, Trans. Am. Math. Soc., 335 (1993), 525-546.

\bibitem{BC} M. Bre\v{s}ar, C. R. Miers, Commutativity preserving mapping of von Neumann algebras, Canad. J. Math., 45 (1993), 659-708.


\bibitem{BM} M. Bre\v{s}ar, C. R. Miers, Strong commutativity preserving maps of semiprime rings, Can. Math. Bull., 37 (1994), 457-460.

\bibitem{DA} Q. Deng, M. Ashraf, On strong commutativity preserving maps. Results Math., 30 (1996), 259-263.

\bibitem{FS} C.-K. Fong, A. R. Sourour, On the operator equation
$\sum_{i=1}^n A_iXB_i\equiv 0$, Canad. J. Math., 31 (1979), 845-857.

\bibitem{L} L. Liu, Strong commutativity preserving maps on von Neumann algebras, Linear and Multilinear Algebra., 63 (2015), 490-496.


\bibitem{LL} J.-S. Lin , C.-K. Liu, Strong commutativity preserving maps on Lie ideals, Linear Algebra Appl., 428 (2008), 1601-1609.

\bibitem{MS} L. Moln\'{a}r, P. \u{S}emrl, Non-linear commutativity preserving maps on self-adjont operators, Q. J. Math., 56(2005), 589-595.

\bibitem{Q} X.- F. Qi, Strong 2-commutativity preserving maps on prime rings.

\bibitem{QH} X.-F. Qi, J.-C. Hou, Nonlinear strong commutativity preserving maps on triangular algebras, Oper. Matrices., 6 (2012), 147-158.

\bibitem{QXH} X.-F. Qi, J.-C. Hou, Nonlinear strong commutativity preserving maps Comm. Algebra., 38 ( 2010), 2790-2796.

\bibitem{Z} J.-H. Zhang, Nonlinear maps preserving Lie products on factor von Neumann algebras, Linear Algebra Appl., 429 (2008), 18-30.




\end{thebibliography}
\end{document}